\documentclass[useAMS,usenatbib,onecolumn]{biom_arxiv}
\usepackage{amsmath,amsfonts}
\usepackage[english]{babel}
\usepackage{epic,latexsym}
\usepackage[dvips]{graphicx}
\usepackage[latin1]{inputenc}
\usepackage[usenames]{color}
\usepackage{setspace}

\def\us{\char`\_}

\def\argmax{\mathop{\rm arg\,max}}
\newcommand{\e}{\varepsilon}
\newcommand{\N}{\mathbb{N}}

\newcommand{\A}{\mathcal{A}}                         % alphabet

\newcommand{\bX}{\mathbf{X}}
\newcommand{\bY}{\mathbf{Y}}
\newcommand{\Fs}{\mathcal{F}_s}

\newcommand{\pr}{\mathbb{P}_\theta} % proba sous theta
\newcommand{\prp}{\mathbb{P}_{\theta'}}             % proba sous theta
\newcommand{\prpi}{\mathbb{P}_\pi}             % proba sous pi
\newcommand{\pro}{\mathbb{P}_{0}} % proba sous un theta0
\newcommand{\esp}{\mathbb{E}_\theta} % esperance sous theta
\newcommand{\espp}{\mathbb{E}_{\theta'}} % esperance sous theta
\newcommand{\espo}{\mathbb{E}_{0}}          % esperance sous theta0
\newcommand{\E}{\mathcal{E}}                       % espace des trajectoires
\newcommand{\dsum}{\displaystyle\sum}

\newtheorem{prop}{Proposition}

\title[Context dependent statistical alignment]{A context dependent pair
  hidden Markov model for statistical alignment}

\author{Ana Arribas-Gil\emailx{ana.arribas@uc3m.es}\\
Departamento  de  Estadística,   Universidad  Carlos  III  de
         Madrid,  C/ Madrid, 126 - 28903 Getafe,  Spain.	   \and 
	   Catherine Matias\emailx{catherine.matias@genopole.cnrs.fr}\\
	   Laboratoire Statistique et Génome, Université d'Évry Val d'Essonne, 
UMR CNRS 8071, USC INRA, \\
523 pl. des Terrasses de l'Agora, 91000 Évry,  France.
	   }

\begin{document}

%\date{{\it Received }. {\it Revised } .\newline 
%{\it Accepted } .}

\pagerange{\pageref{firstpage}--\pageref{lastpage}} \pubyear{}

\volume{}
\artmonth{}
\doi{}

%  This label and the label ``lastpage'' are used by the \pagerange
%  command above to give the page range for the article

\label{firstpage}

\begin{abstract}
This  article proposes a  novel approach  to statistical  alignment of
nucleotide sequences  by introducing a context  dependent structure on
the  substitution process  in  the underlying  evolutionary model.  We
propose to estimate alignments and context dependent mutation rates relying on the observation of two homologous sequences. The procedure is based on a generalized pair-hidden Markov structure, where conditional on the alignment path, the nucleotide sequences follow a Markov distribution. 
We   use   a   stochastic   approximation   expectation   maximization
(\texttt{saem})  algorithm to give  accurate estimators  of parameters
and  alignments.   We  provide  results  both on  simulated  data  and
vertebrate genomes, which are known to have a high mutation
rate from CG dinucleotide.  
In  particular, we establish  that the method  improves the
accuracy of  the alignment  of a human  pseudogene and  its functional
gene.
\end{abstract}
%
%  Please place your key words in alphabetical order, separated
%  by semicolons, with the first letter of the first word capitalized,
%  and a period at the end of the list.
%


\maketitle

\section{Introduction}
Alignment  of  DNA  sequences  is  the  core  process  in  comparative
genomics. An alignment is a mapping of the nucleotides in one sequence
onto  the  nucleotides in  the  other  sequence which captures two important evolutionary processes: the substitution (of a nucleotide, by another one) and the insertion or deletion (indel, of one or several nucleotides). Thus, in this mapping,  two  \emph{matched} nucleotides are  supposed to derive  from a common ancestor  (they are homologous) and the mapping also allows for
\emph{gaps} which may be introduced into  one or the other sequence. Several
powerful alignment algorithms have been developed to align two or more
sequences, a vast majority of them relying on dynamic programming
methods that optimize an alignment score function. While score-based  alignment procedures appear to be very rapid (an important point regarding to the current sizes of the databases), they have at least two major flaws: i) The choice of the score parameters requires some biological expertise and is non objective. This is particularly problematic since the resulting alignments are not robust to this choice \citep{Dewey_2006}. ii) The highest scoring alignment might not be the most relevant from a biological point of view. One should then examine the $k$-highest scoring alignments (with the additional issue of selecting a value for $k$) but this approach often proves impractical because of the sheer number of suboptimal alignments \citep{Waterman_Eggert}. We refer to Chapters 2 and 4 in \cite{Durbin} for an introduction to sequence alignment.\\

Contrarily to score-based  alignment procedures, statistical (or probabilistic) alignment methods rely on an explicit modeling of the two evolutionary
processes at the core of an alignment, namely the  substitution and
the indel  processes. While  there is  a general
agreement  that  substitution  events  are  satisfactorily  modeled
relying on  continuous time Markov  processes on the  nucleotide state
space, very few indel models  have been proposed in the literature and
none of them has become a standard. The first  evolutionary model for biological sequences dealing with indels was formulated by \cite{TKF}, hereafter the TKF91 model. It provides a basis for  performing alignment within a  statistical framework: the \emph{optimal} alignment is obtained as the evolutionary path that maximizes the likelihood of observing the sequences as being derived from a common ancestor (under this model).  This
has at least two major advantages over the score-based method. 
First, in  this context, the  evolutionary parameters (which  play the
role of an underlying scoring function) are not chosen by the user but rather estimated from the data, in a
maximum  likelihood  framework.  Second,  a posterior  probability  is
obtained on the set of alignments: the analysis does not concentrate on the highest probable path and rather provides sample alignments from this posterior distribution. Moreover, a \emph{confidence measure} may be assigned to each  position in  the  alignment. This  is  particularly relevant  as
alignment  uncertainty is  known to  be  a major  flaw in comparative
genomics \citep{Wong_Suchard}.

Statistical  alignment is performed  under the  convenient pair-hidden
Markov model, hereafter pair-HMM \cite[see][]{AGM06,Durbin}, enabling the use of 
generalizations of the expectation-maximization (\texttt{em}) algorithm \citep{Dempster}. 
Some  drawbacks of  the  preliminary TKF91  proposal  have first  been
improved by  the same authors in  what is called the  TKF92 version of
the model \citep{TKF2}. Then, the original models have been later refined in many ways, as for instance in \cite{Ana_Metzler,KnudsenMiyamoto,Metzler,Miklos_03,MiklosLunterHolmes,Miklos_toro}. We underline that an important characteristic of all those models is that the indel and the substitution processes are handled independently. An introduction to statistical alignment may be found in \cite{Lunter_etal}.

It is  important to  note that  the original TKF91  model and  all the
later   developments   are   primarily   focused   on   modeling   the
indel   process.   As   for   the   underlying
substitution model,  pair-HMMs procedures mostly rely  on the simplest
one:        the         one-parameter        Jukes-Cantor        model
\citep{Jukes_Cantor}.  Surprisingly,   while  models  of  substitution
processes  have  made  important  progress  towards  a  more  accurate
description   of   reality   in   the  past   years   \cite[see   for
instance][Chapter  13]{Felsenstein_book}, these developments  have had
almost no impact  in the statistical alignment literature.  One of the
reasons  for that  may be  the existence  of important  biases  found in
sequence alignment and resulting from incorrectly inferred indels.  For  instance,  \cite{Lunter_Rocco}  exhibit  three  major
classes  of  such biases:  (1)  gap wander,  resulting from  the
incorrect placement of gaps due to spurious non homologous similarity;
(2) gap attraction, resulting in the joining of two closely positioned
gaps into one  larger gap; and (3) gap  annihilation, resulting in the
deletion of  two indels of equal  size for a  typically more favorable
representation as substitutions \cite[see also][]{Holmes_Durbin98}. 
This might explain a stronger focus on accurate estimation of indels positions. Another reason could be that since evolutionary parameters are to be estimated in statistical alignments procedures, their number should remain quite low.  

On the contrary, score-based alignment methods offer some variety in the form of the 
match scoring functions, which is the counterpart of the substitution
process parameters when relying on statistical alignment. Note that in this context, using different types of indels  scores  is  hardly  possible, since  the  dynamic  programming approach  limits   those  scores  to   affine  (in  the   gap  length)
penalties.  Nonetheless, the match scores may take different forms. One  of  the   most  recent  improvements  in  score-based
alignment procedures concerns the development of context-based scoring
schemes  as  in  \cite{Gambin06,Gambin07}  \cite[see also][as  a  much
earlier attempt]{Huang}.  In these works,  the cost of  a substitution
depends on the surrounding nucleotides. As a consequence, the score of
an  alignment depends  on  the order  of  editing operations.  Another
drawback of this  approach lies in the issue  of computing a $p$-value
associated to  the alignment significance. Indeed, in  order to assess
the statistical significance of an alignment, one requires the knowledge of the tail distribution of the score statistic, under the null hypothesis of non homology between the sequences. While this distribution is not even fully characterized in the simpler case of site independent scoring schemes allowing gaps, it is currently out of reach when introducing dependencies in the substitution cost.

Nonetheless, allowing for context dependence in the modeling of the evolutionary
processes is important as it captures local sequence dependent fluctuations in substitution rates. It is
well-known for instance that, at least in vertebrates, there is a
higher  level  of methylated  cytosine  in  CG dinucleotides  (denoted
CpG).  This  methylated  cytosine  mutates  abnormally  frequently  to
thymine which results in higher
substitution rates form the pair CpG \citep{Bulmer}. 
At the mutational level, context dependent substitution rates can have a significant impact on prediction accuracy \citep{Siepel_Haussler}. 

Note that allowing for context dependent indel rates is also a
challenging   issue.  For   instance,  microsatellite   expansion  and
contraction is known to be a context dependent process \cite[see for
instance some  recent results in][]{microsat} and suggests  the use of
context dependent indel rates. Such  an attempt has been very recently
made by  \citet{Hickey_Blanchette,Hickey_Blanchette_JCB},  relying on  tree
adjoining grammars.
However, introducing locally dependent  indel rates in an evolutionary
process implies that the  substitution and the indel processes may
not anymore  be handled separately. This independence  between the two
processes is  a key assumption  in limiting the computational  cost of
pair-HMMs and removing  this hypothesis  is beyond  the scope  of the
present   work.   Thus,   we    limit   our   approach   to   modeling
context dependent substitution rates.  \\

In this work, we formulate a context dependent pair hidden
Markov model, so as to  perform statistical alignment of two  sequences allowing for
context dependent substitution rates.  We then propose a method to perform parameter maximum likelihood estimation as well as posterior sampling within this model.  Our
procedure         relies         on        stochastic         versions
\citep{Celeux_Diebolt1,Celeux_Diebolt2} of  the \texttt{em} algorithm.
A   similar   method   has   already   been   successfully   used   in
\cite{Ana_Metzler}   in   a   classical   pair-HMM    (with   no
context dependent substitution rates) in which case it gives results
very similar to Gibbs sampling. Note that \cite{Holmes} also
suggests  the use  of  a  stochastic version  of  \texttt{em} for  
estimating  indel rates in the multiple alignment framework, assuming the alignment is unknown. In the same spirit, \cite{Hobolth2} proposes an \texttt{MCMC-EM} algorithm for the estimation of neighbor-dependent substitution rates from a multiple alignment and a phylogenetic tree. The advantages of
stochastic  versions of \texttt{em}  algorithm is  that they  allow to
perform maximum likelihood estimation  in a hidden variables framework
where  the  maximization  is  not  explicit, the  computation  of  the
likelihood is computationally very expensive, and where the use of the
classical  \texttt{em} algorithm  is not  feasible. To  our knowledge,
this work is, with \citet{Hickey_Blanchette}'s approach, one of the
first  attempts to  perform  a contextual  statistical alignment.\\

This article is organized as follows. 
In Section~\ref{sec:model}, we first describe a context dependent pair
hidden Markov model. In this framework, parameter maximum likelihood estimation
combined with hidden states posterior sampling provides (a posterior
probability  on)  contextual  alignments  of two  sequences  which  is
data-driven in the sense that it does not require any parameter choice. 
Section~\ref{sec:consistency}  briefly  discusses  asymptotic  consistency
results   that   can  be   obtained   on   the  parameter   estimation
procedure. Then,  Section~\ref{sec:algo} describes two  versions of the
\texttt{em} algorithm applied in this framework of contextual statistical
alignment:  the stochastic  expectation maximization  (\texttt{sem}) and
the  stochastic approximation  expectation  maximization (\texttt{saem})
\citep{saem}  algorithms.  The  performance  of these  procedures  are
illustrated  on synthetic  data in  Section~\ref{sec:simus} and  a real
data set is handled in Section~\ref{sec:realdata}.

\section{Model}\label{sec:model}
\begin{figure}[!h]
\begin{center}
\setlength{\unitlength}{0.7mm}
\begin{picture}(55,55)
%Axes
\put(6,6){\line(1,0){46}} \put(6,6){\line(0,1){46}}
\multiput(16,6)(10,0){4}{\dashbox{1}(0,46)}
\multiput(6,16)(0,10){4}{\dashbox{1}(46,0)}
%Lettres
\put(9,0){\Large A}\put(19,0){\Large A}\put(29,0){\Large
T}\put(39,0){\Large G} \put(0,9){\Large C}\put(0,19){\Large
T}\put(0,29){\Large G}\put(0,39){\Large G}
%Chemin
\thicklines
\drawline[5](6,6)(16,16)(26,16)(36,26)(46,36)(46,46)
\linethickness{0.35mm}
\put(16,16){\line(1,0){10}}\put(46,36){\line(0,1){10}}
\end{picture}
\end{center}
\caption{Graphical  representation  of  an  alignment  between  two  sequences
  $X=AATG$ and $Y=CTGG$. The displayed alignment is $ 
  \stackrel{A}  {\text{\it  \scriptsize C}}  \stackrel{A}{\text{  \large -  }}
  \stackrel{T}{\text{\it  \scriptsize T}} \stackrel{G}{  \text{\it \scriptsize
      G}}   \stackrel{\text   {  \large   -   }}{\text{\it  \scriptsize   G}}.
  $}\label{fig:pairHMM} 
\end{figure}
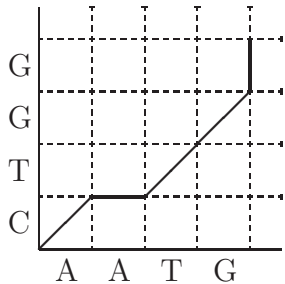

A pair-hidden Markov model is described by the distribution of a non-observed increasing
path through the two-dimensional integer lattice $\N\times \N$ and given this hidden
path,  the conditional distribution of a pair of observed sequences $X_{1:n}:=X_1,\ldots, X_n \in \A^n$ and
$Y_{1:m}:=Y_1,\ldots,Y_m \in  \A^m$, with values on a  finite set $\A$
(the nucleotides alphabet). The path
reflects the  insertion-deletion process  acting on the  sequences and
corresponds to the \emph{bare} alignment, namely the alignment without specification of the nucleotides.
It  is the  sum of  random  moves from  the set  $\mathcal{E}=\{(1,1),
(1,0), (0,1)\}$ also denoted $\{M, I_X,I_Y\}$, where $M$ stands for
'match', $I_X$ stands for 'insert in $X$' and similarly $I_Y$ stands for 'insert in $Y$'. This path represents an empty alignment of the sequences, as the move $M$ represents an homologous site between the two sequences, the move $I_X$ represents either an insertion in the first sequence ($X_{1:n}$) or a deletion in the second ($Y_{1:m}$) and similarly for the move $I_Y$ (see Figure~\ref{fig:pairHMM}). 
 
 We assume that the hidden path follows a Markov distribution without specifying from which insertion-deletion process it comes from. This could be for instance  the  TKF91 insertion-deletion process \citep{TKF} or any of its later variants \citep{Ana_Metzler,KnudsenMiyamoto,Metzler,Miklos_03,TKF2}. 
Given the empty alignment, we then describe the conditional distribution of the sequence letters. In a 'match' position, this distribution reflects the substitution process while in 'insert in X' or 'insert in Y' positions, this distribution reflects both the substitution process and the distribution of inserted segments.\\

More precisely, let $\{\e_t\}_{t\ge 0}$ be an ergodic, homogeneous and
stationary Markov chain with  state space $\E=\{(1,1), (1,0), (0,1)\}$
and transition probability $\pi$ (a stochastic matrix of size $3\times 3$).
%and stationary distribution $\mu$. 
For notational convenience, the
state space $\E$ is sometimes equivalently denoted by $\{M, I_X,I_Y\}$.  
In  this context,  we  might  use obvious  notations  as for  instance
$\pi_{MM}$ for the conditional probability $\mathbb{P}(\e_t=(1,1)|\e_{t-1}=(1,1))$.

Let  $Z_t=(N_t,M_t)=\sum_{s=1}^t \e_s$  be  the path  induced by  this
Markov chain through the lattice $\N\times \N$.
In the following, a path through the lattice $\N \times \N$ is always assumed to be increasing, namely of the form $\sum_s e_s$ for some sequence $e \in \E^{\mathbb{N}}$.
We denote  by $\prpi$ the corresponding stationary  probability and by
$\E_{n,m}$ the set  of paths in $\N\times \N$  starting at $(0,0)$ and
ending at $(n,m)$. For any $e  \in \E_{n,m}$, we also let $|e|$ denote
the length of $e$, which satisfies $n\vee m \le |e|\le n+m$.\\

The hidden process generates the observed sequences in the following way. 
Given   the   hidden   path   $\{\e_t\}_{t\ge   0}$,   the   sequences
$\{ X_{1:n},Y_{1:m}\}$ are generated  according to an order one Markov
chain: for any $e \in \E_{n,m}$,
$$%\forall e \in \E_{n,m}, \quad 
\mathbb{P}(         X_{1:n},Y_{1:m}|          \e         =e)         =
\prod_{s=1}^{|e|}\mathbb{P}(X_{N_s},Y_{M_s}|e_s,e_{s-1},X_{N_{s-1}},Y_{M_{s-1}}).$$
Note that the indexes $N_s,M_s$  of the observations generated at time
$s$ are  random. This is in  sharp contrast with  the classical hidden
Markov model (HMM) and is specific to pair-HMMs \cite[see][for more details]{AGM06}.

In  this   work,  we  restrict  our  attention   to  order-one  Markov
conditional distributions. Note  that the following developments could
be generalized to higher order Markov chains at the cost of increased 
computational burden.

As we want  to model the dependency in  the substitution process only,
we further constrain  the model so that the  dependency occurs only in
successive  homologous (match)  positions. We  thus use  the following
parametrization 
\begin{equation*}
  \mathbb{P}(X_{N_s},Y_{M_s}|e_s,e_{s-1},X_{N_{s-1}},Y_{M_{s-1}}) \\= 
\left\{
  \begin{array}[]{ll}
\tilde{h}(X_{N_s},Y_{M_s}|X_{N_{s-1}},Y_{M_{s-1}})     &\text{if } e_s=e_{s-1}=(1,1), \\
h(X_{N_s},Y_{M_s})     &\text{if } e_s=(1,1)\neq e_{s-1}, \\
f(X_{N_s})     &\text{if } e_s=(1,0), \\
g(Y_{M_s})     &\text{if } e_s=(0,1),
 \end{array}
\right.
\end{equation*}
% \begin{multline*}
%   \mathbb{P}(X_{N_s},Y_{M_s}|e_s,e_{s-1},X_{N_{s-1}},Y_{M_{s-1}}) \\= 
% \left\{
%   \begin{array}[]{ll}
% \tilde{h}(X_{N_s},Y_{M_s}|X_{N_{s-1}},Y_{M_{s-1}})     &\text{if } e_s=e_{s-1}=(1,1), \\
% h(X_{N_s},Y_{M_s})     &\text{if } e_s=(1,1)\neq e_{s-1}, \\
% f(X_{N_s})     &\text{if } e_s=(1,0), \\
% g(Y_{M_s})     &\text{if } e_s=(0,1),
%  \end{array}
% \right.
% \end{multline*}
where $f,g$ are probability measures (p.m.) on $\A$; function $h$ is a p.m. on
$\A^2$ and $\tilde h$ may be viewed as a stochastic matrix on $\A^2$.
%and for any values $a,b\in\A^2$, we have and $\tilde h(\cdot|a,b)$ are p.m. on $\A^2$. 

Obvious  necessary conditions for the parameters  to be identifiable
in this model are the following
\begin{equation}\label{eq:ident}
\left\{
  \begin{array}{cl}
 i)&   \exists   a,b   \in   \A   \text{  such   that   }   h(a,b)\neq
f(a)g(b), \\
ii)& \exists a,b,c,d \in \A \text{ such that }\tilde h(a,b|c,d)\neq
  f(a)g(b), \\
iii)& \exists a,b,c,d \in \A \text{ such that } \tilde h(a,b|c,d)\neq
h(a,b).
  \end{array}
\right .
\end{equation}
The whole parameter set is then given by
\begin{equation*}
\Theta =
 \{\theta =(\pi,f,g,h, \tilde{h}) ; 
\theta \text{ satisfies } \eqref{eq:ident}\} .
\end{equation*}

Statistical alignment  of two  sequences consists in  maximizing with
respect   to  $\theta$   the   following  criterion   $w_{n,m}(\theta)$
\citep{Durbin}. This criterion plays the  role of a log-likelihood  in the
pair-HMM. \cite[For a discussion on the quantities playing the role of
likelihoods in pair-HMMs, see][]{AGM06}.
For  any  integers  $n,m\ge  1$  and any  observations  $X_{1:n}$  and
$Y_{1:m}$, let
\begin{equation}
\label{eq:wnm}
  w_{n,m}(\theta) :=\log Q_\theta(X_{1:n}, Y_{1:m}) \\:=\log \pr (\exists s \geq 1, Z_s =(n,m) ; X_{1:n}, Y_{1:m}).
\end{equation}
% \begin{multline}
% \label{eq:wnm}
%   w_{n,m}(\theta) :=\log Q_\theta(X_{1:n}, Y_{1:m}) \\:=\log \pr (\exists s \geq 1, Z_s =(n,m) ; X_{1:n}, Y_{1:m}).
% \end{multline}
The criterion $w_{n,m}$ is more explicitly defined as 
\begin{multline*}
w_{n,m}(\theta) 
=   \log  \Big(   \sum_{s=n\vee  m}^{n+m}   \sum_{e  \in   \E_{n,m},  |e|=s}
\prpi(\e_{1:s}=e_{1:s}) 
 \times\big\{\prod_{k=1}^s f(X_{n_k})^{1\{e_k=(1,0)\}}
g(Y_{m_k})^{1\{e_k =(0,1)\}} \\
\times h(X_{n_k},Y_{m_k})^{ 1\{e_k=(1,1), e_{k-1}\neq (1,1)\}} 
\tilde{h}(X_{n_k},Y_{m_k}|X_{n_{k-1}},Y_{m_{k-1}})^{ 1\{e_k=e_{k-1}=(1,1)\}}
\big \} \Big) ,
\end{multline*}
% \begin{multline*}
% w_{n,m}(\theta) 
% =   \log  \Big(   \sum_{s=n\vee  m}^{n+m}   \sum_{e  \in   \E_{n,m},  |e|=s}
% \prpi(\e_{1:s}=e_{1:s}) \\
%  \times\big\{\prod_{k=1}^s f(X_{n_k})^{1\{e_k=(1,0)\}}
% g(Y_{m_k})^{1\{e_k =(0,1)\}} \\
% h(X_{n_k},Y_{m_k})^{ 1\{e_k=(1,1), e_{k-1}\neq (1,1)\}} \\
% \tilde{h}(X_{n_k},Y_{m_k}|X_{n_{k-1}},Y_{m_{k-1}})^{ 1\{e_k=e_{k-1}=(1,1)\}}
% \big \} \Big) ,
% \end{multline*}
where $(n_k,m_k) =\sum_{t=1}^k e_t$. 
Now, we define the parameter estimator as 
\begin{equation*}
  \hat \theta_{n,m}=\argmax_{\theta\in \Theta} w_{n,m}(\theta). 
\end{equation*}
In  the  classical  pair-HMM  (namely  without  accounting  for
dependencies in the substitution process), consistency results for the
above    estimator    $\hat    \theta_{n,m}$    were    obtained    in
\cite{AGM06} and  are extended here  to the context dependent  case in
Section~\ref{sec:consistency}.  Then in 
Section~\ref{sec:algo}, we develop an \texttt{saem} algorithm to compute $\hat \theta_{n,m}$ and provide an a posteriori probability
distribution over  the set of alignments.  Indeed, a main  issue is to
obtain an accurate alignment  of the sequences. Score-based alignment
methods  heavily   rely  on  the   Viterbi   algorithm  for  this
purpose \cite[see][]{Durbin}, which provides the most probable \emph{a posteriori} alignment. However,
optimal  alignments  often  look  different  from  typical  ones,  and
providing a unique alignment as the result of the estimation procedure
may not be very informative. For  this reason, it may be more interesting to provide the probability
distribution of hidden states over each pair of nucleotides from the observed
sequences (for the parameter value $\hat \theta_{n,m}$). This can be done under the pair-HMM framework and gives us a reliability measure of any alignment of the two sequences \cite[see][]{Ana_Metzler}.

\section{Consistency results for parameter estimation}\label{sec:consistency}

In this  part, we  give generalizations of  results first  obtained in
\cite{AGM06}, to the context dependent pair-HMM. As the proofs follow the
same  lines  as in  this  reference, we  shall  omit  them.  We  first
introduce some  notations and formulate some assumptions  that will be
needed. For any $\delta >0$, we let
\begin{equation*}
  \Theta_\delta =\{\theta \in \Theta ; \forall k,
  \theta_k \ge \delta\} \quad \text{ and } \quad\Theta_0=\{\theta \in
  \Theta ; \forall k, \theta_k >0 \} .
\end{equation*}
The true parameter value will  be denoted by $\theta_0$ and is assumed
to belong to the set $\Theta_{0}$. Probabilities and expectations under this parameter value are denoted
by $\pro$ and $\espo$, respectively. 
We now introduce a notation  for the marginal distributions of $h$ and
$\tilde h$. For any $b,c\in \A$, let
\begin{align*}
&h_{X}(\cdot) :=\sum_{a\in \A}h(\cdot,a) , \;
h_{Y}(\cdot) :=\sum_{a\in \A}h(a,\cdot) ,\\
&\tilde h_{X}(\cdot|b,c) :=\sum_{a\in \A}\tilde h(\cdot,a|b,c) ,\;
\tilde h_{Y}(\cdot|b,c) :=\sum_{a\in \A} \tilde h(a,\cdot|b,c) .
\end{align*}
We then let $\Theta_{marg}$ be a subset of parameters in $\Theta_0$ satisfying
some assumptions on the marginals of $h$ and $\tilde h$. 
\begin{equation*}
  \Theta_{marg}=\Big\{\theta\in\Theta_0; \forall b,c \in \A, 
  h_{X} = f,   h_{Y} =g, \\
  \tilde h_{X}(\cdot|b,c) =f(\cdot), 
  \tilde h_{Y}(\cdot|b,c)=g(\cdot)
\Big\}.
\end{equation*}
Note that for instance, the parametrization used in Section~\ref{sec:simus}
satisfies  $\theta  \in  \Theta_{marg}$.   Finally, we  introduce  the
process $w_t$ defined in a similar way as in~\eqref{eq:wnm} by 
$$
w_t(\theta):=\log Q_\theta (X_{1:N_t},Y_{1:M_t}).
$$
As previously noted, the length  $t$ of an alignment of sequences with
respective sizes $n$ and $m$ satisfies $n\vee
m\le  t\le n+m$.  Thus, asymptotic  results for  $t \to  +\infty$ will
imply  equivalent   ones  for   $n,m\to  +\infty$.  In   other  words,
consistency results obtained  when $t \to +\infty$ are  valid for long
enough observed sequences, even if one does not know the length $t$ of the true underlying alignment. We now state the main theorem. 
\begin{theorem}
\label{thm:consistency}
 For any $\theta\in\Theta_0$, we have
\begin{itemize}
\item [i)] Renormalized log-likelihood $ t^{-1} w_{t} (\theta)$ converges $\pro$-almost surely and in
$\mathbb{L}_1$, as $t$ tends to infinity, to
\begin{equation*}
 w (\theta ) = \lim_{t\rightarrow \infty}\frac{1}{t}\espo \left(\log
Q_\theta(X_{1:N_{t}},Y_{1:M_{t}}) \right) \\= \sup_{t}\frac{1}{t}\espo
\left(\log Q_\theta(X_{1:N_{t}},Y_{1:M_{t}} ) \right).
\end{equation*}
\item  [ii)]  Moreover,  $w  (\theta_0)\ge  w (\theta)  $  and  strict
  inequality  is valid  as  long as at  least  one of  these
  conditions is satisfied 
  \begin{itemize}
  \item [A.] $\theta \in \Theta_{0}$ and $\forall \lambda >0, \esp (\e_1) \neq \lambda \espo(\e_1)$,
\item [B.] $\theta_0, \theta\in \Theta_{marg}$ and either $f\neq f_0$ or $g\neq g_0$.
  \end{itemize}
\item [iii)] The  family of functions $\{t^{-1}w_t(\theta)\}_{t\ge 1}$
  is  uniformly equicontinuous  on  the set  $\Theta_\delta$, for  any
  $\delta >0$.
\end{itemize}
\end{theorem}
The  three  statements in  Theorem~\ref{thm:consistency}  are the  key
ingredients  to  establish  the  consistency of  $\hat  \theta_{n,m}$,
following  Wald's  classical proof \citep{Wald} of  maximum likelihood  estimators
consistency \cite[see also][Chapter 5]{VDV}.
The most subtle part is point $ii)$, that states some sufficient conditions under
which the maximum of the limiting function $w$ is attained only at the
true parameter  value $\theta_0$.  For instance, condition  $A$ means
that  the main directions  of the  paths $\{Z_t\}_{t\ge  1}$ generated
under the two parameter  values $\theta$ and $\theta_0$ are different.
As  a consequence, these  two parameters  may be  easily distinguished
from each  other relying on  the distributions $\pr$ and  $\pro$. Note
that these two different cases (condition $A$ or $B$ satisfied) are certainly
not  exhaustive and we  refer to  \cite{AGM06} for  simulation results
indicating that point $ii)$ might be true under other scenarios. 
Note  in particular  that for evolutionary models assuming that
insertions and  deletions happen at the  same rate, $\esp(\epsilon_1)$
is proportional to $(1,1)$ for any parameter value $\theta$. In this case, condition $A)$ is never satisfied and we are not able to prove that $h\neq
h_0$  or $\tilde  h  \neq  \tilde h_0$  implies  strict inequality  $w
(\theta_0) > w (\theta) $. This  is due to the dependency structure on
the  sequences, that  make  it  difficult to  link  the difference  $w
(\theta_0)- w  (\theta) $  with a Kullback-Leibler  divergence between
densities $h$ and $h_0$ or between $\tilde h$ and $\tilde h_0$ \cite[see the proof of Theorem 2 in][for more details]{AGM06}.\\

In general, one is not interested in the full (context) pair-HMM
parameters but rather in  a sub-parametrization induced by the choice
of   a   specific   evolutionary   model   (see   for   instance   the
sub-parametrization proposed  in Section~\ref{sec:simus}). Let $\beta
\mapsto \theta(\beta)$ be a  continuous parametrization from some set
$B$     to     $\Theta$.      For     any     $\delta     >0$,     let
$B_\delta=\theta^{-1}(\Theta_\delta)$.       We       assume      that
$\beta_0:=\theta^{-1}(\theta_0)$   belongs  to  $B_\delta$   for  some
$\delta >0$. We then let 
$$
\hat \beta_{n,m} :=\argmax_{\beta \in B_{\delta}} w_{n,m}(\theta(\beta)).
$$
Besides  from  giving  results  on  the  frequentist  estimator  $\hat
\beta_{n,m}$, we consider Bayesian  estimates of the parameters. Let $\nu$
be a  prior probability  measure on the  set $B_\delta$.  Markov chain
Monte  Carlo  (MCMC) algorithms  approximate  the random  distribution
$\nu_{n,m}$, interpreted  as the posterior  measure given  observations $X_{1:n}$
and $Y_{1:m}$ and defined by 
\begin{equation*}
\nu_{n,m} (d\beta) = 
\frac{Q_{\theta(\beta)} (X_{1:n},Y_{1:m})\nu(d\beta)}
{\int_{B_{\delta}} Q_{\theta(\beta ')} (X_{1:n},Y_{1:m})\nu(d\beta')}. 
\end{equation*}
With these definitions at hand, we can now state the following corollary.
\begin{corollary}
  If  the maximizer of  $\beta \to w(\theta(\beta))$ over  the set
    $B_\delta$ is  reduced to the singleton  $\{\beta_0\}$, then we have the following results
  \begin{itemize}
  \item [i)] $\hat \beta_{n,m}$ converges $\pro$-almost-surely to $\beta_0$,
\item [ii)] If $\nu$ puts some weight on a neighborhood of $\beta_0$, then
  the sequence of posterior measures $\nu_{n,m}$ converges in distribution, $\pro$-almost surely, to the Dirac mass at $\beta_0$.
  \end{itemize}
\end{corollary}

\section{Algorithms}\label{sec:algo}
The approach  we use here to maximize   criterion~\eqref{eq:wnm} and provide a posterior distribution on the set of alignments, given   two  observed   sequences   relies  on   a  stochastic   version \citep{Celeux_Diebolt1,Celeux_Diebolt2}  of  \texttt{em}  algorithm
\citep{BPSW,Dempster}. This approach has already been successfully used
in \cite{Ana_Metzler} in a classical pair-HMM (meaning with no
context dependent substitution rates) in which case it gives results
very  similar  to  Gibbs   sampling.  It  allows  to  perform  maximum
likelihood estimation in a framework  in which the maximization is not
explicit, the  computation of  the likelihood is  computationally very
expensive,  and in  which  the  use of  \texttt{em}  algorithm is  not
feasible (see Appendix for details).

\subsection{Forward and backward equations}
We  describe below the  forward and
backward  probabilities  used in  the later  procedures.   The formulas  are
easily obtained by generalizing those from the classical pair-HMM \citep{Ana_thesis,Durbin}.

The forward  equations are obtained  as follows. For any  value $u\in
\E$,  let $\E_{n,m}^{u}=\{e\in  \E_{n,m},  e_{|e|}=u\}$ be  the set  of
paths ending  at $(n,m)$  with last step  being $u$. Then  the forward
probabilities are defined for any $i\le n$ and $j\le m$, as
\begin{equation*}
\alpha^u(i,j)=\pr(\exists s\ge 1, Z_s=(i,j), \e_s=u, X_{1:i},Y_{1:j}) \\= \sum_{e\in \E_{n,m}^u}\pr(\e=e,X_{1:i},Y_{1:j} ).
\end{equation*}

These forward probabilities are computed recursively in the following way. We first set the initial values
\begin{eqnarray*}
 & \forall u\in \E,\quad
\alpha^u(0,-1)=\alpha^u(-1,0)=0
\\
 & \forall u\in \{I_X,I_Y\}, \quad \alpha^u(0,0)=0 \text{ and } \alpha^{M}(0,0)=1.
\end{eqnarray*}
Then for $i=0,\ldots, n$ and $j=0,\ldots,m$ except $(i,j)=(0,0)$, we recursively compute
\begin{multline}
  \label{eq:forward_match}
 \alpha^{M}(i,j) =\pr(\exists s\ge 1, \e_s=(1,1), Z_s=(i,j), X_{1:i},Y_{1:j})  
=\sum_{u\in \E} \pr (\exists s\ge 1, \e_s=(1,1), Z_s=(i,j), \e_{s-1}=u, X_{1:i},Y_{1:j}) \\
=  \tilde{h}(X_i,Y_j|X_{i-1},Y_{j-1})\pi_{MM}\alpha^{M}(i-1,j-1) 
+h(X_i,Y_j)\pi_{I_XM}\alpha^{I_X}(i-1,j-1) 
+h(X_i,Y_j)\pi_{I_YM}\alpha^{I_Y}(i-1,j-1),
\end{multline}
and in the same way
\begin{eqnarray*}
  \alpha^{I_X}(i,j) &=& f(X_i) \sum_{u\in \E} \pi_{uI_X} \alpha^u(i-1,j) , \\%\quad
  \alpha^{I_Y}(i,j) &=& g(Y_j) \sum_{u\in \E} \pi_{uI_Y} \alpha^u(i,j-1) .
\end{eqnarray*}
Note  that the  log-likelihood of  the observed  sequences is  then simply
obtained from the forward probabilities by taking
\begin{equation*}
  w_{n,m}(\theta)=\log \pr(\exists s\ge 1, Z_s=(n,m), X_{1:n},Y_{1:m}) \\
=\log (\alpha^{I_X}(n,m)+\alpha^{I_Y}(n,m)+\alpha^{M}(n,m)).
\end{equation*}

Note  that the forward  equations could  thus be  used to  compute the
log-likelihood of  two observed sequences and a  numerical optimization of
this quantity could give maximum likelihood estimates. However, such a
strategy fails short  as soon as the dimension  of the parameter space
is  not small.   Indeed,  the computational  cost  of calculating  the
log-likelihood  for a  given  parameter  value is  $O(nm)$,  and thus  the
computational  cost of  its maximization  over a  grid  is $O(p^knm)$,
where $k$ is the number of parameters and $p$ the number of values considered for each parameter.  \\

We now describe the backward probabilities and their recursive computation. For any value $u\in \E$, the backward probabilities are defined as
\begin{equation*}
  \beta^u(i,j)=  \pr(X_{i+1:n},Y_{j+1:m}|\exists  s\ge  1,  Z_s=(i,j),
  \e_s=u, X_i,Y_j).
\end{equation*}
We initialize these values by taking, for any $u \in \E$, 
\begin{eqnarray*}
\beta^u(n,m)  =  1 \text{  and }  \beta^u
  (n,m+1)=\beta^u(n+1,m)= 0,
\end{eqnarray*}
and recursively compute, for $i=n,\ldots, 1$ and $j=m,\ldots,1$, except $(n,m)$, the quantities
\begin{multline}\label{eq:backward}
    \beta^u(i,j) =\pi_{uM} \beta^{M}(i+1,j+1) 
%\\\times 
\{ h(X_{i+1},Y_{j+1})1_{u\in\{I_X,I_Y\}} +\tilde h(X_{i+1},Y_{j+1}|X_i,Y_j) 1_{u=M} \} \\
+\pi_{uI_X}\beta^{I_X}(i+1,j)       f(X_{i+1})       +      \pi_{uI_Y}
\beta^{I_Y}(i,j+1) g(Y_{j+1}) .
\end{multline}

Note that contrarily to the case of HMM, the forward-backward
probabilities do not give rise to the classical  \texttt{em} strategy, as we do not obtain from these equations the conditional expectation of the complete log-likelihood, given the two observed sequences  (see Appendix for further details).

\subsection{\texttt{sem} and \texttt{saem} algorithms}\label{algos} 
The forward-backward probabilities may  be used to simulate the hidden path,
conditional on the observed sequences. Thus, \texttt{em} algorithm may
be replaced by \texttt{sem} \citep{Celeux_Diebolt1} or \texttt{saem} \citep{saem}
procedures. Indeed,  let us explain  the basic idea of  these stochastic
versions    of   \texttt{em}    algorithm.
Denoting by $L_{n m}$  the random value $s\ge 1$  such that $Z_s=(n,m)$
(the first and  only hitting time for the point  $(n,m)$, which is not
necessarily finite), the complete log-likelihood is 
$$\log \pr(X_{1:n}, Y_{1:m}, L_{n m},\e_{1:L_{n m}}) $$
(see Appendix for details) and for a current parameter value $\theta'$, its conditional
expectation writes 
$$Q_{\theta'}(\theta):=\espp(\log \pr(X_{1:n}, Y_{1:m}, L_{n m},\e_{1:L_{n m}}) |X_{1:n}, Y_{1:m}).$$ 
Now, the  idea of  these stochastic approximations  is to  replace the
computation of $Q_{\theta'}(\theta)$ by an approximation obtained from
the simulation (under parameter value $\theta'$) of a number of hidden paths $\{\e_{1:s}\}$ in
$\E_{n,m}$, with possibly different lengths. 
These  paths satisfy  the property  that  for each  obtained value  of
$s$, the paths with length $s$ have the same occurrence probability as $\prp(\e | X_{1:n}, Y_{1:m}, L_{n m}=s)$. Note that the length of each simulated path is random. 
Then the maximization with respect to $\theta$ of the complete log-likelihood is done on the basis of the simulated alignments. As a consequence, this maximization is performed via a simple counting of specific events and does not require numerical procedures. \\
In particular, iteration $r$ of  \texttt{sem} algorithm writes
\begin{itemize}
\item[-] {\em Simulation step:} generate one hidden path $e^r$ of some random length $s$ with the same distribution as $${\mathbb P}_{\theta^{(r-1)}}(\e_{1:s} | X_{1:n}, Y_{1:m}, L_{n m}=s).$$
\item[-]  {\em Maximization  step:}  $\theta^{(r)}=\argmax_{\theta \in
    \Theta} \log \pr(X_{1:n},$ $ Y_{1:m}, $ $\e=e^r)$.
\end{itemize}
And iteration $r$ of  \texttt{saem} algorithm writes
\begin{itemize}
\item[-] {\em Simulation step:} generate $m(r)$ hidden paths $e^r (j)$, $j=1,\dots,m(r)$ with resulting random lengths $s_{r,j}$, each one having the same distribution as $${\mathbb P}_{\theta^{(r-1)}}(\e_{1:s_{r,j}} | X_{1:n}, Y_{1:m}, L_{n m}=s_{r,j}).$$
\item[-] {\em Stochastic approximation step:} update 
  \begin{equation*}
 Q_r       (\theta)      =      Q_{r-1}       (\theta)     \\ +\gamma_r
\Big(\dfrac{1}{m(r)}\sum_{j=1}^{m(r)}   \log   \pr(X_{1:n},  Y_{1:m},
  \e=e^r(j)) 
- Q_{r-1}(\theta) \Big),
  \end{equation*}
where $\{\gamma_r\}_{r\geq 1}$ is a decreasing sequence of positive step size and 
$Q_0(\theta)= 1/m(0)\sum_{j=1}^{m(0)} $ $\log \pr(X_{1:n}, Y_{1:m}, \e=e^0(j))$ .
\item[-] {\em Maximization step:} $\theta^{(r)}=\argmax_{\theta \in \Theta} Q_r (\theta) $.
\end{itemize}

The  simulation step  is common  to  both algorithms  and consists  in
drawing paths  of some  unknown length $s$  with same  distribution as
${\mathbb   P}_{\theta^{(r-1)}}(\e_{1:s}  |  X_{1:n},   Y_{1:m},  L_{n
  m}=s)$. This  is performed via  the backwards sampling based  on the
forward probabilities as we now explain. We may write 
\begin{eqnarray}
{\mathbb P}_{\theta^{(r-1)}}(\e_{1:s} | X_{1:n}, Y_{1:m}, L_{n m}=s) 
&= &{\mathbb P}_{\theta^{(r-1)}}(\e_s| X_{1:n},Y_{1:m},L_{n m}=s) 
\prod_{k=1}^{s-1} {\mathbb P}_{\theta^{(r-1)}}(\e_k | \e_{k+1:s}, X_{1:n},Y_{1:m},L_{n m}=s)  \nonumber\\
&=&{\mathbb P}_{\theta^{(r-1)}}(\e_s| X_{1:n},Y_{1:m},Z_s=(n,m)) \prod_{k=1}^{s-1} {\mathbb P}_{\theta^{(r-1)}}(\e_k | \e_{k+1}, Z_{k+1}, X_{1:N_{k+1}}, Y_{1:M_{k+1}}) .\label{eq:backwards_sampling} 
\end{eqnarray}
% \begin{eqnarray}
% &&{\mathbb P}_{\theta^{(r-1)}}(\e_{1:s} | X_{1:n}, Y_{1:m}, L_{n m}=s)  \nonumber\\
% &=& {\mathbb P}_{\theta^{(r-1)}}(\e_s| X_{1:n},Y_{1:m},L_{n m}=s) \nonumber\\
% &&\times \prod_{k=1}^{s-1} {\mathbb P}_{\theta^{(r-1)}}(\e_k | \e_{k+1:s}, X_{1:n},Y_{1:m},L_{n m}=s)  \nonumber\\
% &=&{\mathbb P}_{\theta^{(r-1)}}(\e_s| X_{1:n},Y_{1:m},Z_s=(n,m)) \nonumber\\
% &&\times \prod_{k=1}^{s-1} {\mathbb P}_{\theta^{(r-1)}}(\e_k | \e_{k+1}, Z_{k+1}, X_{1:N_{k+1}}, Y_{1:M_{k+1}}) .\label{eq:backwards_sampling} 
% \end{eqnarray}
The last equality comes from the following: Conditional on $\e_{k+1}$,
the  random   variable  $\e_{k}$  does  not  anymore   depend  on  the
observations  after  time  $k+2$.  If  we moreover  condition  on  the
variables  $ \e_{k+1:s}$  and  $L_{n m}=s$, then  we  know which  point
$(i,j)$        on       the        lattice        corresponds       to
$Z_{k+1}=(n,m)-\sum_{t=k+2}^s\e_t=(i,j)$.  The values of  the observed
sequence   after  time   $k+2$  correspond   exactly   to  $X_{i+1:n},
Y_{j+1:m}$.       Then       conditional       on       $\{\e_{k+1:s},
X_{1:n},Y_{1:m},L_{n m}=s\}$,  the variable  $\e_{k}$  only depends  on
$\e_{k+1}$ (according to the Markov property) and on the values $Z_{k+1}, X_{1:N_{k+1}}, Y_{1:M_{k+1}}$.

So, given $\e_{k+1}$, the position $Z_{k+1}=(i,j)$ on the lattice and the observed values $X_{1:i}, Y_{1:j}$, we want to sample $\e_k$ from ${\mathbb P}_{\theta^{(r-1)}}(\e_k | \e_{k+1}, X_{1:i}, Y_{1:j}, Z_{k+1}=(i,j))$, for $k=s,s-1,\ldots,1$. Up to a renormalizing constant, each step corresponds to sampling from ${\mathbb P}_{\theta^{(r-1)}}(\e_k, \e_{k+1} |X_{1:i},$ $ Y_{1:j},$ $ Z_{k+1}=(i,j))$ or equivalently from $\pi_{\e_k,\e_{k+1}}\alpha^{\e_k}((i,j)-\e_{k+1})$. 
Then the backwards sampling writes 
\begin{itemize}
\item[-] $e^r_s$ is sampled to be $v$ with probability 
$$
\dfrac{\alpha^v(n,m)}{\alpha^M(n,m)+\alpha^{I_X}(n,m)+\alpha^{I_Y}(n,m)} ,$$
\item[-] For $k<s$ and $Z_{k+1}=(i,j)$, $i>1, j> 1$, if $e^r_{k+1}=M$,
  then $e^r_k$ is sampled to be $M$ with probability 
$$\dfrac{\alpha^M(i-1,j-1)\pi_{MM}
  \tilde{h}(X_{i},Y_{j}|X_{i-1},Y_{j-1}) }{\alpha^M(i,j)} ,$$
and to be $u=I_X$ or $I_Y$ with probability 
$$\dfrac{\alpha^u(i-1,j-1)\pi_{uM} h(X_{i},Y_{j}) }{\alpha^M(i,j)} .$$
Note that these latter probabilities are correctly renormalized according to Equation~\eqref{eq:forward_match}.
\item[-] For $k<s$ and $Z_{k+1}=(i,j)$, $i>1, j> 1$,  if $e^r_{k+1}=I_X$ then $e^r_k$ is sampled to be $u$ with probability 
  \begin{equation}
    \label{eq:sample_IX}
\dfrac{\alpha^u(i-1,j)\pi_{uI_X} f(X_{i}) }{\alpha^{I_X}(i,j)} ,     
  \end{equation}
\item[-] For $k<s$ and $Z_{k+1}=(i,j)$, $i>1, j> 1$,  if $e^r_{k+1}=I_Y$ then $e^r_k$ is sampled to be $u$ with probability 
\begin{equation}
  \label{eq:sample_IY}
\dfrac{\alpha^u(i,j-1)\pi_{uI_Y} g(Y_{j}) }{\alpha^{I_Y}(i,j)} ,   
\end{equation}
\item[-] For $k<s$ and $Z_{k+1}=(1,j)$, $ j>1$, if $e^r_{k+1}=M$ or $I_X$ then all the values $e^r_l, 1\le l\le k$ are fixed to $I_Y$, otherwise $e^r_k$  is sampled to be $u$ with probability given by \eqref{eq:sample_IY}.
\item[-] For $k<s$ and $Z_{k+1}=(i,1)$, $ i>1$, if $e^r_{k+1}=M$ or $I_Y$ then all the values $e^r_l, 1\le l\le k$ are fixed to $I_X$, otherwise $e^r_k$  is sampled to be $u$ with probability given by \eqref{eq:sample_IX}.
\end{itemize}
According to Equation~\eqref{eq:backwards_sampling}, we thus have sampled paths $e^r$ satisfying the following property.
\begin{prop}
  The paths $e \in \E_{n,m}$ sampled from the above scheme with parameter value $\theta$ occur with probability 
$\pr(\e=e|X_{1:n},Y_{1:m},L_{n m}=|e| )$.
\end{prop}

We now consider the complete log-likelihood obtained from the observed
sequences and one simulated path. We have 
\begin{multline*} 
\log  \pr(X_{1:n},   Y_{1:m},\e =e^r)  =   \sum_{u\in \E} 1_{e^r_1=u}\log  \pi_{u}  +\sum_{k=2}^{|e^r|}\sum_{u,v\in     \E^2} 1_{e^r_{k-1}=u,e^r_{k}=v} \log\pi_{uv}\\
 +\sum_{k=1}^{|e^r|}\sum_{a\in \A} [1_{e^r_k=I_X, X_{n_k}=a} \log  f(a) +1_{e^r_k =I_Y,Y_{m_k=a}} \log g(a)]
+ \sum_{k=1}^{|e^r|}  \sum_{a,b\in \A^2} 1_{e^r_k=M,  X_{n_k}=a, Y_{m_k}=b}
\log h(a,b)\\
+  \sum_{k=2}^{|e^r|}   \sum_{a,b,c,d\in  \A^4}  1_{e^r_k=e^r_{k-1}=M,
  (X_{n_k}, Y_{m_k},X_{n_{k-1}}, Y_{m_{k-1}}) =(a,b,c,d)} \log \tilde h(a,b|c,d) ,
\end{multline*}
% \begin{eqnarray*} 
% &&\log \pr(X_{1:n}, Y_{1:m},\e  =e^r) = \sum_{u\in \E} 1_{e^r_1=u}\log
% \pi_{u} \\ &&+\sum_{k=2}^{|e^r|}\sum_{u,v\in \E^2} 1_{e^r_{k-1}=u,e^r_{k}=v} \log\pi_{uv}\\
%  &&+\sum_{k=1}^{|e^r|}\sum_{a\in \A} [1_{e^r_k=I_X, X_{n_k}=a} \log  f(a) +1_{e^r_k =I_Y,Y_{m_k=a}} \log g(a)]\\
% &&+ \sum_{k=1}^{|e^r|}  \sum_{a,b\in \A^2} 1_{e^r_k=M,  X_{n_k}=a, Y_{m_k}=b}
% \log h(a,b)\\
% &&+  \sum_{k=2}^{|e^r|}   \sum_{a,b,c,d\in  \A^4}  1_{e^r_k=e^r_{k-1}=M,
%   (X_{n_k}, Y_{m_k},X_{n_{k-1}}, Y_{m_{k-1}}) =(a,b,c,d)}\\
% &&\times \log \tilde h(a,b|c,d) ,
% \end{eqnarray*}
where we recall that $(n_k,m_k)= \sum_{i=1}^k e_i^r$ (for simplicity, we omit the subscript $r$).
Then  in the  maximization  step of  \texttt{sem},  the components  of
$\theta$ are  updated to the values:  $ \forall u,v  \in \E$, $a,b,c,d
\in \A $,
\begin{align}
&   \pi^{(r)}  _{u}=\dfrac{N_{e^r}(u)}{\dsum_{v\in   \E}   N_{e^r}(v)},   \quad
\pi^{(r)}_{u,v}=\dfrac{N_{e^r}(u,v)}{\dsum_{w\in   \E}   N_{e^r}(u,w)}, \quad
f^{(r)}(a)=\dfrac{N_{e^r}(a|I_X)}{\dsum_{a'\in\A}     N_{e^r}(a'|I_X,)}, \quad 
 g^{(r)}(a)=\dfrac{N_{e^r}(a|I_Y)}{\dsum_{a'\in\A}
  N_{e^r}(a'|I_Y,)},      \nonumber \\
&h^{(r)}(a,b)=\dfrac{N_{e^r}(a,b|M)}{\dsum_{a',b'\in                 \A}
  N_{e^r}(a',b'|M)}, 
 \quad \tilde{h}^{(r)}(a,b|c,d)=\dfrac{N_{e^r}(a,b|MM;c,d)}{\dsum_{a',b'\in \A}
  N_{e^r}(a',b'|MM;c,d)}, 
\label{Mstep}
\end{align}
% \begin{align}
% &   \pi^{(r)}  _{u}=\dfrac{N_{e^r}(u)}{\dsum_{v\in   \E}   N_{e^r}(v)},   \quad
% \pi^{(r)}_{u,v}=\dfrac{N_{e^r}(u,v)}{\dsum_{w\in   \E}   N_{e^r}(u,w)}, \nonumber \\ 
% &f^{(r)}(a)=\dfrac{N_{e^r}(a|I_X)}{\dsum_{a'\in\A}     N_{e^r}(a'|I_X,)}, \quad 
%  g^{(r)}(a)=\dfrac{N_{e^r}(a|I_Y)}{\dsum_{a'\in\A}
%   N_{e^r}(a'|I_Y,)},      \nonumber \\
% &h^{(r)}(a,b)=\dfrac{N_{e^r}(a,b|M)}{\dsum_{a',b'\in                 \A}
%   N_{e^r}(a',b'|M)}, \nonumber \\
% % \quad
% &\tilde{h}^{(r)}(a,b|c,d)=\dfrac{N_{e^r}(a,b|MM;c,d)}{\dsum_{a',b'\in                 \A}
%   N_{e^r}(a',b'|MM;c,d)}, 
% \label{Mstep}
% \end{align}
where  we use  obvious  notations  for the  counts  $N_e$ in  sequence
$e$. Namely for any $ u,v \in \E$, and any $a,b,c,d \in\A$, we let 
 $N_{e}(u)= \sum_k  1\{e_k=u\}$, $N_{e}(u,v)=\sum 1\{e_{k-1}=u,e_k=v\}$,
 $N_{e}(a|I_X)=\sum   1\{e_k=I_X,X_{N_k}=a\}$    and   similarly   for
 $N_{e}(a|I_Y)$, $N_{e}(a,b|M)$ and also $ N_{e}(a,b|MM; c,d)=\sum 1\{e_k=e_{k-1}=M,(X_{N_k},Y_{M_k},X_{N_{k-1}},Y_{M_{k-1}}) =(a,b,c,d)\}$.

In the maximization step of \texttt{saem} we have to take into account the $m(r)$ hidden paths and the stochastic approximation of $Q_r (\theta)$. Then the values of $\theta$ are updated as in~\eqref{Mstep} by replacing the counts $N_{e^r}(\cdot)$ by their counterparts $\tilde{N}_r(\cdot)$, where $\tilde{N}_r=\tilde{N}_{r-1} +\gamma_r \left(\dfrac{1}{m(r)}\sum_{j=1}^{m(r)} N_{e^r(j)} - \tilde{N}_{r-1}\right)$.

From a theoretical point of view, the convergence properties of the two algorithms are different. In \texttt{saem}, the sequence $\{\theta^r\}_{r\ge 1}$ converges, under general conditions, to a local maximum of the log-likelihood. It is important to note that the stochastic perturbation introduced with respect to the original \texttt{em} algorithm allows avoiding saddle points, which are possible attractive stationary points of the sequence generated by \texttt{em} \citep{saem}. In \texttt{sem}, the sequence $\{\theta^r\}_{r\ge 1}$ does not converge pointwise, but it is an homogeneous Markov chain which is ergodic under appropriate conditions, see for instance \cite{Diebolt96}.

In  this  context,  maximum  likelihood estimators  are  sensitive  to
overfitting if  there is  insufficient data. Indeed,  if an  event has
never been  observed in the sequences, the  corresponding estimator is
not well defined, as both the numerator and denominator in~\eqref{Mstep} are zero. 
Pseudo-counts \cite[Chapter  3]{Durbin} or other  smoothing techniques
such as  the one  proposed in \cite{Kneser_Ney}  may be used  to solve
this problem.  However, when using \texttt{saem}  algorithm to perform
maximum  likelihood  estimation,  this   problem  is  minimized  and  no
smoothing strategies  are required. Indeed, since at  iteration $r$ we
generate $m(r)$ hidden paths, the counts $\tilde{N_r}$ (computed as an
average over  those alignments) are  rarely equal to $0$,  even during
the   first  iterations   when   $\gamma_r$  is   typically  set   to
$1$. However, if we want  to completely avoid this possibility, we may
replace  event  counts by  pseudo-counts  in  the  case of  unobserved
events, during these  first iterations. Note also that  in practice we
do not generate whole alignments at each iteration step, but only 
sub-alignments   of  a   certain  given   length  within   the  current
alignment.   This   strategy  combines   an   MCMC  procedure   within
\texttt{saem} and  is known to preserve   convergence properties of
  \texttt{saem}  algorithm    \citep[see][for  further
details]{Ana_Metzler}. In terms of event counts, this strategy implies that a current low
estimate of an event occurrence probability will have some
impact on next iteration estimate. Nonetheless, at the last iterations
of the  algorithm, when typically  $\gamma_r < 1$, current counts
are averaged with counts from previous iterations, which provides a smoothing procedure by itself.

\section{Simulations}\label{sec:simus}

For these simulations, we consider nucleotide sequences with alphabet 
$\A=\{A,C,G,T\} $  and a constrained  parameter set, motivated  by the
choice of an underlying evolutionary model on these sequences. 
As for  the generation of  the hidden
Markov chain, we  consider the TKF91 indel model  with equal value for
insertion   and   deletion    rate   $\lambda$   \cite[see][for   more
details]{Ana_thesis}. The transition matrix of the hidden path is thus given by 
\begin{equation*} 
\pi = \dfrac{1}{1+\lambda}  
\begin{pmatrix}
 e^{-\lambda} & 1-e^{-\lambda} & \lambda\\ 
\dfrac{\lambda e^{-\lambda}}{1-e^{-\lambda}} & \lambda& 1+\lambda-\dfrac{\lambda}{1-e^{-\lambda}} \\ 
e^{-\lambda} & 1-e^{-\lambda}& \lambda
\end{pmatrix} ,
\end{equation*}
where the order of the states in the matrix is $M$, $I_X$ and $I_Y$.
Then, we describe the conditional distribution of the sequences. For any $x\in \A$, we set $f(x)=g(x)=\mu_x$, where $\mu$ is going to be the stationary distribution of the substitution process, according to the parametrization of $h$ and $\tilde{h}$ considered below.
We set the distribution $h$ according to a simple substitution process
with  equal  substitution  rate  $\gamma$  and  different  nucleotides
frequencies (modified Jukes Cantor). Let 
\begin{equation*}
  h(x,y) =
\left\{
  \begin{array}{ll}
\mu_x(1-e^{-\gamma})\mu_y & \text{if } x\neq y\\
\mu_x [\mu_x(1-e^{-\gamma})+e^{-\gamma} ] &\text{if } x=y.   
  \end{array}
\right.
\end{equation*}
The  parameter   $\tilde  h$  accounting   for  the  context   in  the
substitution process differs from $h$ only in a $\substack{ C \\ C}$ match context. In this way, we want to take into account a possibly higher substitution rate from the pair CpG to CpA. More precisely, we let
\begin{equation*}
 \mathbb{P}(X_{N_s}=x,Y_{M_s}=y|e_s=e_{s-1}=M,\\
X_{N_{s-1}}=x',Y_{M_{s-1}}=y')= \left\{
     \begin{array}{ll}
   h_C(x,y) & \text{if } x'=y'=C,\\
  h(x,y) &  \text{otherwise} .      
     \end{array}
\right.
\end{equation*}
In a $\substack{ C \\ C}$ match context, we use a model initially introduced in two slightly different forms by \cite{Felsenstein84} and \cite{HKY}. This model considers different substitution rates for transitions (\emph{i.e.} substitutions within the chemical class of purines $\mathcal{R}=\{A,G\}$ or pyrimidines $\mathcal{Y}=\{C,T\}$) and for transversions (substitutions modifying the chemical class) as well as different nucleotide frequencies $\mu$. 
Denoting by $\bar  x$ the other nucleotide in  the same chemical class
as  $x$  (namely   $x\neq  \bar  x$  and  either   $\{x,\bar  x\}  \in
\mathcal{R}$ or $\{x,\bar x\} \in \mathcal{Y}$), we get 
\begin{equation*}
   h_C(x,y)= \left\{
    \begin{array}{ll}
\mu_x e^{-(\alpha +\beta)} +\mu_xe^{-\beta} (1-e^{-\alpha})\frac{\mu_x}{\mu_x+\mu_{\bar x}} +\mu_ x(1-e^{-\beta})\mu_x & \text{ if } y=x ,\\
\mu_xe^{-\beta}(1-e^{-\alpha})\frac{\mu_{\bar x}}{\mu_{\bar x}+\mu_{x}}
+\mu_x(1-e^{-\beta})\mu_{\bar x}& \text{ if } y=\bar x, \\
\mu_x(1-e^{-\beta})\mu_y & \text{ otherwise. } 
    \end{array}
\right.
\end{equation*}
% \begin{align*}
%  &h_C(x,x)=\mu_x  e^{-(\alpha  +\beta)}  +\mu_xe^{-\beta}
%   (1-e^{-\alpha})\frac{\mu_x}{\mu_x+\mu_{\bar x}} \\
% &\qquad +\mu_x(1-e^{-\beta})\mu_x ,\\
% &h_C(x,\bar x)= \mu_xe^{-\beta}(1-e^{-\alpha})\frac{\mu_{\bar x}}{\mu_{\bar x}+\mu_{x}}
% +\mu_x(1-e^{-\beta})\mu_{\bar x} ,\\
% &\text{Otherwise, }  h_C(x,y)=\mu_x(1-e^{-\beta})\mu_y .
% \end{align*}
We refer to \cite{Felsenstein_book} for more details on substitution processes.\\

\begin{figure}[h!]
\begin{center}
\includegraphics[width=12cm,height=3cm]
%[trim = 25mm 120mm 25mm 110mm, clip,width=12cm]
{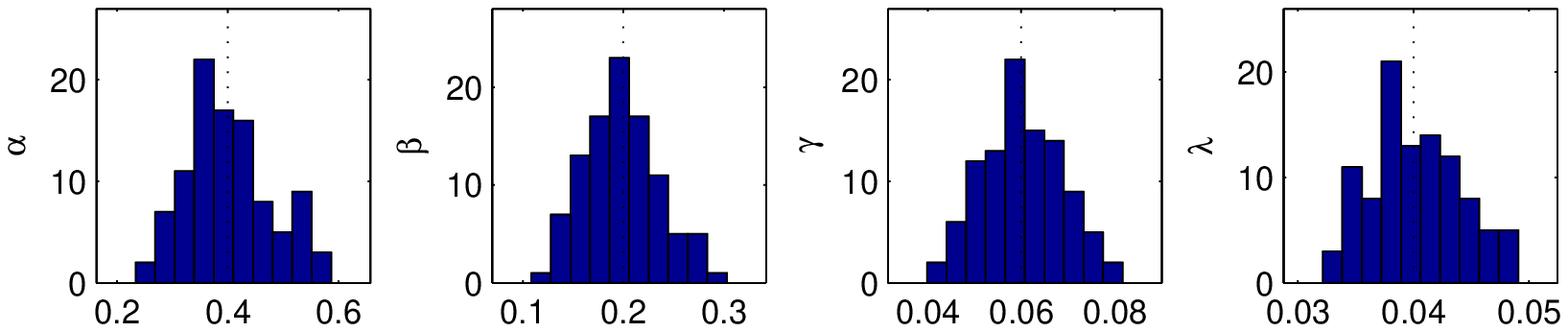}\\
\includegraphics[width=12cm,height=3cm]
%[trim = 25mm 120mm 25mm 120mm, clip,width=12cm]
{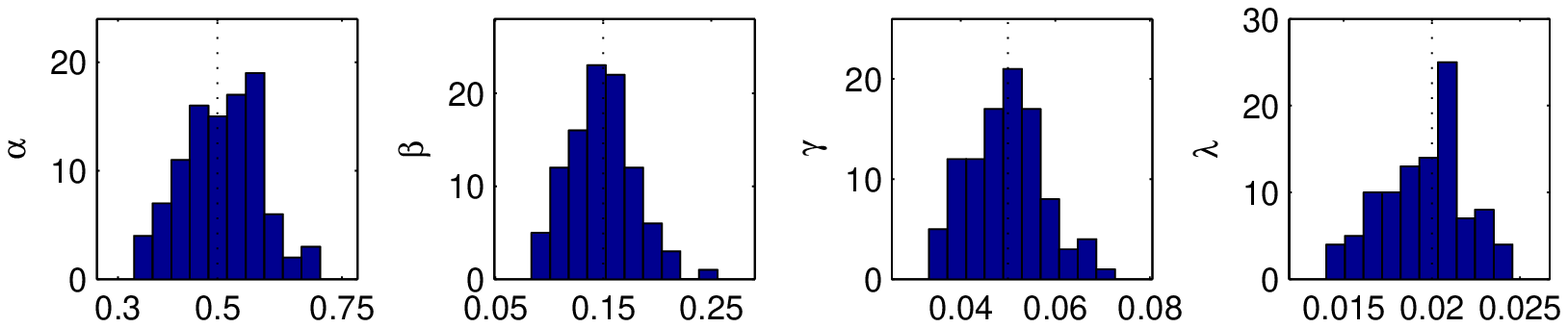}\\
\caption{Distribution  of  \texttt{saem}   estimates  of  the  reduced
  parameter  vector  $(\alpha,  \beta, \gamma, \lambda)$  over  100
  simulation runs for the two simulated  data sets (Data set 1 on top,
  Data set 2 on bottom). The real values of
  the parameters are displayed in dotted line.}
\label{F1}
\end{center}
\end{figure}

\begin{figure}[h!]
\begin{center}
\includegraphics[width=12cm,height=3cm]{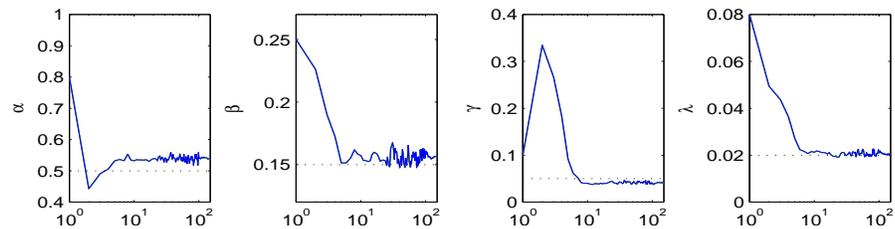}
%[trim = 25mm 110mm 25mm 110mm, clip,width=12cm]{taux_param_log.eps}
\caption{\texttt{saem}  estimates  of  the  reduced  parameter  vector
  $(\alpha, \beta, \gamma,  \lambda)$ over iterations in a single
  simulation run from the second data set. A logarithmic scale is used
  for the $x$-axis.  The true parameter values are displayed in
  dotted line.}
\label{F2}
\end{center}
\end{figure}

According to this parametrization we  have simulated two sets of data
with an alignment length of 2000 base pairs (bp) and parameter values described in
Table~\ref{tab:models}.
\begin{table}
\caption{Parameter values in the two data sets.}  
\label{tab:models}
\begin{tabular}{c | cccc}
\\
Data set 1 & $\alpha=0.4$ &$\beta=0.2$&$\gamma=0.06$ &$\lambda=0.04$ \\
\hline \\
Data set 2 &$\alpha=0.5$& $\beta=0.15$ & $\gamma=0.05$ &$\lambda=0.02$ \\
 \end{tabular}
\bigskip 
\end{table}
In these data sets, the  substitution rate $\gamma$ is larger than the
indel  rate $\lambda$  as expected  in biological  sequences,  and the
transition and  transversion substitution rates in a  $\substack{ C \\
  C}$ match context ($\alpha$ and $\beta$) are larger than the regular
substitution rate $\gamma$. In both sets of data we have set the stationary distribution of the substitution process to
\begin{equation*}
\mu_A=0.225 \quad \mu_C=0.275 \quad \mu_G=0.275 \quad \mu_T=0.225 ,
\end{equation*}
in order to  obtain a GC content (proportion of Gs or Cs) larger than
$50\%$.\\ 
For  each  one of  the  data sets  we  have  conducted two  estimation
procedures. On  the one  hand, we have  estimated the  whole parameter
$(\pi,  f,  g,  h,  \tilde{h})$  via the  \texttt{saem}  algorithm  as
explained in Section~\ref{algos}. On  the other hand, we have combined
the  \texttt{saem} algorithm with  numerical optimization  to estimate
only  the evolutionary  parameters $\lambda$,  $\gamma$,  $\alpha$ and
$\beta$.  The  first procedure  has  the  advantages  of being  robust
against  misspecification  of the  underlying  evolutionary model  and
relying  on  a explicit  maximization  step  based  on the  counts  of
events.  The  second  procedure  is  computationally  more  expensive,
however,  it is more  parsimonious and  it provides  a straightforward
evolutionary interpretation of the  parameter values. In both cases we
have performed 150 iterations  of the \texttt{saem} algorithm with the
parameter $\gamma_r$  of the stochastic  approximation set to  $1$ for
$r=1,\dots, 100$ and to  $1/(r-100)$ for $r=101,\dots,150$. The number
of simulated hidden paths is $m(r)=5$, $r=1,\dots, 20$ and $m(r)=10$,
$r=21,\dots, 150$. We  used the same initial values  of the parameters
for  every simulation  run  over  the two  simulated  data sets.  These
initial  values,  estimates  and  standard  deviations  are  given  in
Tables~\ref{T1a},~\ref{T1b} (first procedure) and~\ref{T2} (second procedure). In
Figure~\ref{F1} we present the distribution of estimates obtained with
the second  estimation approach  for the two  simulated data  sets. We can observe that both procedures provide unbiased and very accurate estimates. In
Figure~\ref{F2} we show the convergence of the \texttt{saem} estimates
on  a  single simulation  run  from  the  second data  set.  We use a logarithmic scale for the $x$-axis to highlight the very fast convergence of the parameters to their true values. 
\begin{table}
\caption{Mean  values  and  standard  deviations  of  estimates  over  100
  simulation runs for Data Set 1. The estimation has been
  performed  without  taking into
  account the specific parametrization used to generate the data.}
\label{T1a}
\begin{center}
\begin{tabular}{|l|c|c|c|}
          \hline
           & Initial value & True value & Estimate  (sd)     \\ 
 \hline 
 $\pi_{MM}$  & 0.85          & 0.9238     & 0.9184 (0.0101)   \\
 $\pi_{MI_X}$& 0.075         & 0.0377     & 0.0407 (0.0074)  \\
 $\pi_{MI_Y}$& 0.075         & 0.0385     & 0.0410 (0.0077) \\
 $\pi_{I_XM}$& 0.85          & 0.9424     & 0.8807 (0.0998)   \\
 $\pi_{I_XI_X}$& 0.075       & 0.0385     & 0.0387 (0.0308)  \\
 $\pi_{I_XI_Y}$& 0.075       & 0.0191     & 0.0805 (0.1034) \\
 $\pi_{I_YM}$& 0.85          & 0.9238     & 0.8788 (0.0982)  \\
 $\pi_{I_YI_X}$& 0.075       & 0.0377     & 0.0821 (0.0987)  \\
 $\pi_{I_YI_Y}$& 0.075       & 0.0385     & 0.0391 (0.0305)  \\
 \hline
 $h(A,A)$      & 0.0625      & 0.2148     & 0.2155 (0.0109)  \\
 $h(A,C)$      & 0.0625      & 0.0036     & 0.0031 (0.0019) \\
 $h(A,G)$      & 0.0625      & 0.0036     & 0.0032 (0.0019)  \\
 $h(A,T)$      & 0.0625      & 0.0029     & 0.0027 (0.0018)  \\
 $h(C,A)$      & 0.0625      & 0.0036     & 0.0034 (0.0021)  \\
 $h(C,C)$      & 0.0625      & 0.2634     & 0.2628 (0.0130)   \\
 $h(C,G)$      & 0.0625      & 0.0044     & 0.0044 (0.0024) \\
 $h(C,T)$      & 0.0625      & 0.0036     & 0.0034 (0.0021)  \\
 $h(G,A)$      & 0.0625      & 0.0036     & 0.0030 (0.0017)   \\
 $h(G,C)$      & 0.0625      & 0.0044     & 0.0041 (0.0023)   \\
 $h(G,G)$      & 0.0625      & 0.2634     & 0.2653 (0.0120)  \\
 $h(G,T)$      & 0.0625      & 0.0036     & 0.0032 (0.0018)   \\ 
 $h(T,A)$      & 0.0625      & 0.0029     & 0.0026 (0.0016)   \\
 $h(T,C)$      & 0.0625      & 0.0036     & 0.0031 (0.0017)   \\
 $h(T,G)$      & 0.0625      & 0.0036     & 0.0035 (0.0021)   \\
 $h(T,T)$      & 0.0625      & 0.2148     & 0.2166 (0.0113)   \\ 
 \hline
 $h_C(A,A)$    & 0.0625      & 0.1600     & 0.1626 (0.0169)  \\
 $h_C(A,C)$    & 0.0625      & 0.0112     & 0.0108 (0.0053)   \\
 $h_C(A,G)$    & 0.0625      & 0.0446     & 0.0454 (0.0111)   \\
 $h_C(A,T)$    & 0.0625      & 0.0092     & 0.0091 (0.0051)   \\
 $h_C(C,A)$    & 0.0625      & 0.0112     & 0.0102 (0.0056)  \\
 $h_C(C,C)$    & 0.0625      & 0.2055     & 0.2072 (0.0193)   \\
 $h_C(C,G)$    & 0.0625      & 0.0137     & 0.0126 (0.0056)   \\
 $h_C(C,T)$    & 0.0625      & 0.0446     & 0.0429 (0.0108)  \\
 $h_C(G,A)$    & 0.0625      & 0.0446     & 0.0440 (0.0098)   \\
 $h_C(G,C)$    & 0.0625      & 0.0137     & 0.0140 (0.0061)   \\
 $h_C(G,G)$    & 0.0625      & 0.2055     & 0.2050 (0.0174)   \\
 $h_C(G,T)$    & 0.0625      & 0.0112     & 0.0110 (0.0053)   \\ 
 $h_C(T,A)$    & 0.0625      & 0.0092     & 0.0086 (0.0047)  \\
 $h_C(T,C)$    & 0.0625      & 0.0446     & 0.0456 (0.0106)   \\
 $h_C(T,G)$    & 0.0625      & 0.0112     & 0.0111 (0.0062)   \\
 $h_C(T,T)$    & 0.0625      & 0.1600     & 0.1598 (0.0197)   \\ 
\hline \\
\bigskip 
\end{tabular} 
\end{center}
\end{table}

\begin{table}
\caption{Mean  values  and  standard  deviations  of  estimates  over  100
  simulation runs for Data set 2. The estimation has been
  performed  without  taking into  account the specific parametrization used to generate the data.}
\label{T1b}
\begin{center}
\begin{tabular}{|l|c|c|c|}
       \hline
           & Initial value & True value & Estimate  (sd)        \\
 \hline 
 $\pi_{MM}$  & 0.85        & 0.9610     & 0.9559 (0.0071) \\
 $\pi_{MI_X}$& 0.075       & 0.0194     & 0.0221 (0.0051) \\
 $\pi_{MI_Y}$& 0.075       & 0.0196     & 0.0219 (0.0055) \\
 $\pi_{I_XM}$& 0.85         & 0.9706     & 0.8477 (0.1262) \\
 $\pi_{I_XI_X}$& 0.075      & 0.0196     & 0.0199 (0.0221) \\
 $\pi_{I_XI_Y}$& 0.075      & 0.0098     & 0.1323 (0.1289) \\
 $\pi_{I_YM}$& 0.85         & 0.9610     & 0.8407 (0.1321) \\
 $\pi_{I_YI_X}$& 0.075      & 0.0194     & 0.1367 (0.1379) \\
 $\pi_{I_YI_Y}$& 0.075 & 0.0196 & 0.0226 (0.0252) \\ 
 \hline
 $h(A,A)$      & 0.0625    & 0.2165     & 0.2168 (0.0114) \\
 $h(A,C)$      & 0.0625    & 0.0030     & 0.0027 (0.0015) \\
 $h(A,G)$      & 0.0625    & 0.0030     & 0.0027 (0.0014) \\
 $h(A,T)$      & 0.0625    & 0.0025     & 0.0024 (0.0015) \\
 $h(C,A)$      & 0.0625    & 0.0030     & 0.0026 (0.0017) \\
 $h(C,C)$      & 0.0625    & 0.2653     & 0.2673 (0.0122) \\
 $h(C,G)$      & 0.0625    & 0.0037     & 0.0033 (0.0020) \\
 $h(C,T)$      & 0.0625    & 0.0030     & 0.0027 (0.0016) \\
 $h(G,A)$      & 0.0625    & 0.0030     & 0.0023 (0.0015) \\
 $h(G,C)$      & 0.0625    & 0.0037     & 0.0033 (0.0020) \\
 $h(G,G)$      & 0.0625    & 0.2653     & 0.2643 (0.0096) \\
 $h(G,T)$      & 0.0625     & 0.0030     & 0.0023 (0.0013) \\ 
 $h(T,A)$      & 0.0625     & 0.0025     & 0.0024 (0.0015) \\
 $h(T,C)$      & 0.0625     & 0.0030     & 0.0025 (0.0016) \\
 $h(T,G)$      & 0.0625     & 0.0030     & 0.0027 (0.0015) \\
 $h(T,T)$      & 0.0625     & 0.2165     & 0.2196 (0.0102) \\ 
 \hline
 $h_C(A,A)$    & 0.0625      & 0.1588     & 0.1595 (0.0159) \\
 $h_C(A,C)$    & 0.0625      & 0.0086     & 0.0069 (0.0038) \\
 $h_C(A,G)$    & 0.0625      & 0.0505     & 0.0504 (0.0096) \\
 $h_C(A,T)$    & 0.0625      & 0.0071     & 0.0067 (0.0042) \\
 $h_C(C,A)$    & 0.0625      & 0.0086     & 0.0087 (0.0045) \\
 $h_C(C,C)$    & 0.0625      & 0.2053     & 0.2053 (0.0199) \\
 $h_C(C,G)$    & 0.0625      & 0.0105     & 0.0103 (0.0046) \\
 $h_C(C,T)$    & 0.0625      & 0.0505     & 0.0504 (0.0101) \\
 $h_C(G,A)$    & 0.0625      & 0.0505     & 0.0515 (0.0102) \\
 $h_C(G,C)$    & 0.0625      & 0.0105     & 0.0097 (0.0052) \\
 $h_C(G,G)$    & 0.0625      & 0.2053     & 0.2036 (0.0191) \\
 $h_C(G,T)$    & 0.0625      & 0.0086     & 0.0078 (0.0043) \\ 
 $h_C(T,A)$    & 0.0625      & 0.0071     & 0.0078 (0.0041) \\
 $h_C(T,C)$    & 0.0625      & 0.0505     & 0.0515 (0.0098) \\
 $h_C(T,G)$    & 0.0625      & 0.0086     & 0.0084 (0.0040) \\
 $h_C(T,T)$    & 0.0625      & 0.1588     & 0.1613 (0.0166) \\ 
\hline \\
\bigskip 
\end{tabular} 
\end{center}
\end{table}

\begin{table}
\caption{Mean  values  and  standard  deviations  of  estimates  over  100
  simulation  runs for  the two  data  sets. The  estimation has  been
  performed by numerical optimization  on the reduced parameter vector
  $(\alpha, \beta,\gamma,\lambda)$.}
\label{T2}
\begin{center}
\begin{tabular}{|l|c|c|c|}
           \multicolumn{4}{|c|}{Data set 1} \\
          \hline
           & Initial value & True value & Estimate   (sd)     \\
 \hline 
 $\alpha$    & 0.8         & 0.4        & 0.400 (0.0750)   \\
 $\beta $    & 0.25        & 0.2        & 0.1985 (0.0372)   \\
 $\gamma$    & 0.1         & 0.06       & 0.0602 (0.0087)   \\
 $\lambda$   & 0.08        & 0.04       & 0.0402 (0.0040)  \\
\hline\\
%%%%%
 \multicolumn{4}{|c|}{Data set 2} \\
          \hline
           & Initial value & True value & Estimate   (sd)     \\
 \hline 
 $\alpha$    & 0.8         & 0.5      & 0.5091 (0.0772) \\
 $\beta $ & 0.25 & 0.15 & 0.1493 (0.0296) \\
 $\gamma$    & 0.1      & 0.05     & 0.0499 (0.0079) \\
 $\lambda$   & 0.08        & 0.02     & 0.0196 (0.0024) \\
\hline\\
\bigskip 
\end{tabular} \end{center}
\end{table}


\section{Application to real data}\label{sec:realdata}
\begin{figure}[h!]
\begin{center}
\includegraphics[width=12cm,height=7cm]
%[trim = 20mm 100mm 20mm 80mm, clip,width=12cm]
{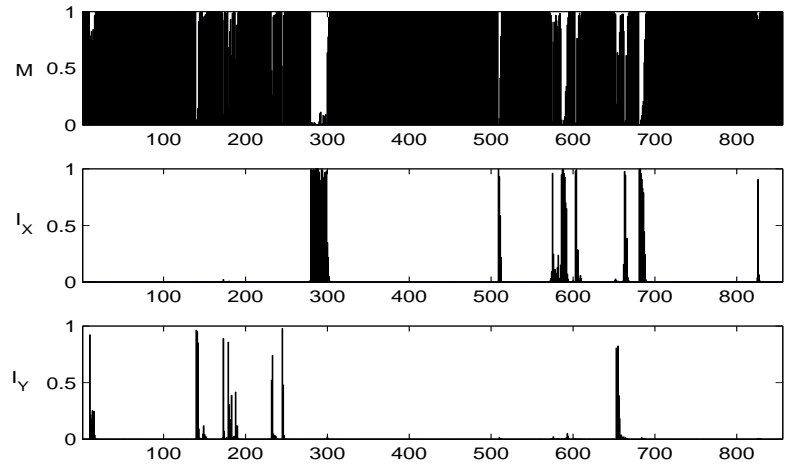}
\end{center}
\caption{Posterior probability distribution of the alignment states at the maximum likelihood parameter estimate with the substitution context-dependent model.}\label{aligCC}
\end{figure}

In this section we illustrate our method through an alignment of the human alpha-globin pseudogene (HBPA1) with its functional counterpart, the human alpha-globin gene (HBA1). This example is inspired by the work of \cite{Bulmer}, who studied the neighbouring base effects on substitution rates in a series of vertebrate pseudogenes, including the human alpha-globin, concluding that there is an increase in the frequency of substitutions from the dinucleotide CpG.\\
We have  extracted the sequences of  the gene and  the pseudogene, who
are located in  the human alpha-globin gene cluster  on chromosome 16,
from   the  UCSC  Genome   Browser  database   \citep{GenomeBrowser}.
 We  have  considered the  whole  sequences
(coding and  non-coding sequence in  the case of the  functional gene)
for the alignment, obtaining sequences  lengths of 842 bp for HBA1 and
812 bp for HBPA1. Because of  the presence of introns and exons in the
HBA1  sequence, it is  natural to  consider a  model allowing  for two
different substitution  behaviors along the  sequence. That is  why we
have used  a model for which  the state space of  the hidden Markov
chain is  $\{M_1,M_2,I_X,I_Y\}$, where  $M_i$, $i=1,2$, stand  for two
different  match   states,  with  different   associated  substitution
processes  $h_i, \tilde{h}_i$, $i=1,2$.  This kind  of model  may also
handle  fragments  insertion  and   deletion  in  contrast  to  single
nucleotide  indels \citep[see][for  more details]{Ana_Metzler}.  As in
the  simulation studies of  the previous  section, $\tilde{h}$  is the
substitution matrix in  a $\substack{ C \\ C}$  match context, and $h$
is the substitution matrix in any  other case. In this way, we want to
take into account a possibly  higher transition rate from nucleotide G
occurring  in the  dinucleotide CpG.  In order  to avoid  any possible
misspecification of  the underlying indel  and substitution processes,
we  have  decided not  to  parametrize  any  of those,  conducting  an
estimation approach equivalent to the first procedure described in the
previous section.  The only difference  here lies in the  dimension of
the    model    parameter,     since    this    parameter    is    now
$(\pi,f,g,h_1,h_2,\tilde{h}_1,\tilde{h}_2)$,   where    $\pi$   is   a
$4\times4$  transition   probability  matrix,  and  $h_i,\tilde{h}_i$,
$i=1,2$ are four different $4\times4$ stochastic vectors, yielding a total number of $12+3+3+4\times 15=78$ free parameters.\\
We have run  \texttt{saem} algorithm on the sequences, performing 500 iterations with $\gamma_r$ set to $1$ for
$r=1,\dots, 400$ and to $1/(r-400)$ for $r=401,\dots,500$.  The number
of simulated hidden paths is $m(r)=5$, $r=1,\dots, 20$ and $m(r)=10$,
$r=21,\dots, 500$. In  Table~\ref{transrate}, we present the estimated
substitution probabilities from nucleotide G. As expected, we found an
increase in  the substitution occurrence for  G in a  $\substack{ C \\
  C}$ match context. We present in Figure~\ref{aligCC} the posterior probability distribution over
the set of possible alignment columns, match, insertion and deletion (the two match states have been merged) for every pair of aligned positions in a consensus alignment obtained at the convergence stage (last iterations) of  \texttt{saem} algorithm.

\begin{table}
\caption{Relative substitution probabilities for G, namely $\phi_i(G,\cdot)/\sum_{a\in\A} \phi_i(G,a)$, $i=1,2$, $\phi=h,\tilde{h}$. The probability of transition for G increases from $0.08$ in the general context to $0.88$ in the $ ^{C}_{ C}$ match context in the first substitution regime, and from $0.02$ to $0.15$ in the second.}
\label{transrate}
\begin{tabular}{ccccc}
\\
\hline
           & A & C & G & T \\
\hline
$h_1$    &  0.0839  &0.0759&   0.8398 &   0.0004\\ 
$\tilde{h}_1$      &   0.8835 &   0.0093&   0.0951&   0.0121\\
$h_2$      &   0.0242&  0.0596&   0.8681&  0.0481\\
$\tilde{h}_2$      &  0.1508&   0.1420&   0.7058&   0.0014 \\
\hline\\
\bigskip 
\end{tabular}
\end{table}

\begin{figure}
\begin{center}
\includegraphics[width=12cm,height=7cm]
%[trim = 20mm 100mm 20mm 80mm, clip,width=12cm]
{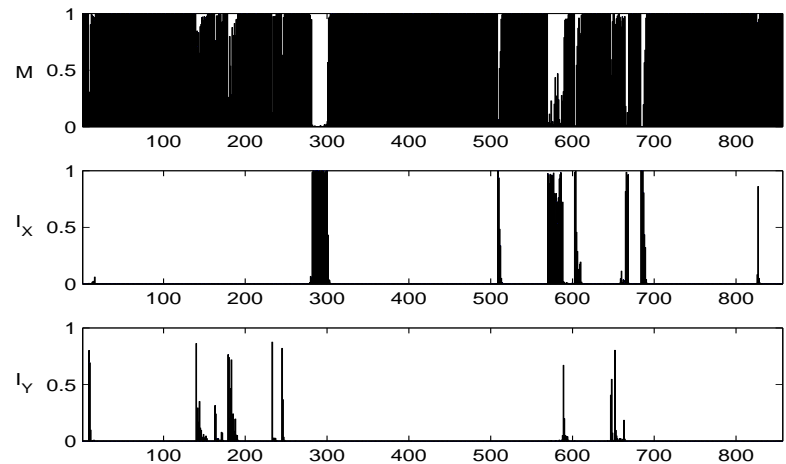}
\end{center}
\caption{Posterior probability distribution of the alignment states at the maximum likelihood parameter estimate without taking into account possible substitution context-dependence.}
\label{alig}
\end{figure}

In order to compare our method to a traditional pair-HMM alignment, we have also run  \texttt{saem} algorithm, with the same settings as before, on a simpler model without substitution context-dependence. In this model we consider again two different substitution regimes, so the parameter vector is now $(\pi,f,g,h_1,h_2)$, yielding a total number of $12+3+3+2\times 15=48$ free parameters. We present in Figure~\ref{alig} the posterior probability distribution over
the set of possible alignment columns, match, insertion or deletion (the two match states have been merged) for every pair of aligned positions in a consensus alignment obtained at the convergence stage (last iterations) of  \texttt{saem} algorithm. This posterior probability is quite similar to the one obtained with the context-dependent model.\\
In  order to  further compare  the two  models, we  have  calculated a
BIC-type criterion penalizing the  maximum likelihood by the number of
parameters.  Specifically, we have  computed $BIC=-2  L +  k \log(n)$,
where $L$  is the value of the  log-likelihood at $\hat{\theta}_{ML}$,
$n$  is  the  sample  size  and  $k$  is  the  total  number  of  free
parameters. We have taken $n=842$, the length of the longest sequence,
since the  length of  the alignment is  unknown, but is  equivalent to
those of  the observed sequences. Indexing by  1 the context-dependent
model and by 2 the non context-dependent one, we get
\begin{eqnarray*} 
 BIC_{1}&=& -2\cdot (-7128.50)+ 78 \log(842) = 14782.39        \\
 BIC_{2}&=& -2\cdot (-7337.40)+ 48 \log(842) = 14998.12,  
\end{eqnarray*}
which means that the context-dependent model better fits this data set.

\section{Discussion}
We have proposed a model for statistical sequence alignment accounting
for context-dependence in the substitution process. This model extends
the pair hidden Markov model in the sense that, conditional on the alignment, the presence of particular nucleotides at the different sites of a sequence are no longer independent. We have adapted the traditional dynamic programing algorithms to the new framework, and we have proposed an efficient estimation method based on  \texttt{saem} algorithm. We have obtained asymptotic results for maximum likelihood estimation on the new model, and through a simulation study, we have shown the accuracy of our algorithms on finite sample estimation. Finally, we have illustrated our methods with the alignment of the human alpha-globin pseudogene and its functional counterpart. We have compared the new model with a classical pair-HMM through a model selection criterion, concluding that taking into account the context-dependence on the substitution process improves the fitting of the data.

\appendix

\section{}
\subsection{The need for a stochastic approximation}
Let us explain why  \texttt{em} algorithm may not be applied in pair-HMM. Denoting  by $L_{n m}$ the random  value $s\ge 1$  such  that $Z_s=(n,m)$  (the  first and  only hitting time for the point $(n,m)$, which is not necessarily finite), the complete log-likelihood writes
\begin{align*}
&\log \pr(X_{1:n}, Y_{1:m}, L_{n m},\e_{1:L_{n m}}) =  \sum_{s=n\vee
  m}^{n+m} 1_{L_{n m}=s} 
 \Big\{ \sum_{u\in \E}
1_{\e_1=u}\log \pi_{u} +\sum_{k=2}^{s}\sum_{u,v\in     \E^2}
1_{\e_{k-1}=u,\e_{k}=v} \log\pi_{uv} \\
& +\sum_{k=1}^s\sum_{a\in \A} [
1_{\e_k=I_X, X_{N_k}=a} \log  f(a) +1_{\e_k =I_Y,Y_{M_k=a}} \log g(a)]
+ \sum_{k=1}^s  \sum_{a,b\in \A^2} 1_{\e_k=M,\e_{k-1}\neq M,  X_{N_k}=a, Y_{M_k}=b} \log h(a,b)
\\  &  +  \sum_{k=2}^s  \sum_{a,b,c,d\in  \A^4}  1_{\e_k=\e_{k-1}=M  ,
  (X_{N_k}, Y_{M_k} , X_{N_{k-1}}, Y_{M_{k-1}})=(a,b,c,d) }  \log \tilde h(a,b|c,d) \Big\}.
\end{align*}
% \begin{align*}
% &\log \pr(X_{1:n}, Y_{1:m}, L_{n m},\e_{1:L_{n m}}) =  \sum_{s=n\vee
%   m}^{n+m} 1_{L_{n m}=s} \\
% & \times \Big\{ \sum_{u\in \E}
% 1_{\e_1=u}\log \pi_{u} +\sum_{k=2}^{s}\sum_{u,v\in     \E^2}
% 1_{\e_{k-1}=u,\e_{k}=v} \log\pi_{uv} \\
% & +\sum_{k=1}^s\sum_{a\in \A} [
% 1_{\e_k=I_X, X_{N_k}=a} \log  f(a) +1_{\e_k =I_Y,Y_{M_k=a}} \log g(a)]
% \\
% &+ \sum_{k=1}^s  \sum_{a,b\in \A^2} 1_{\e_k=M,\e_{k-1}\neq M,  X_{N_k}=a, Y_{M_k}=b} \log h(a,b)
% \\  &  +  \sum_{k=2}^s  \sum_{a,b,c,d\in  \A^4}  1_{\e_k=\e_{k-1}=M  ,
%   (X_{N_k}, Y_{M_k} , X_{N_{k-1}}, Y_{M_{k-1}})=(a,b,c,d) } \\
% & \qquad \times \log \tilde h(a,b|c,d) \Big\}.
% \end{align*}
To     simplify     notations,     we    let     $\bX:=X_{1:n},
\bY:=Y_{1:m}, L:=L_{n m}$ and $\Fs=\{\bX,\bY, L=s\}$. Moreover, we let $(\e,X,Y)_k:=(\e_k,X_{N_k},Y_{M_k})$.
Taking  the  expectation  of  the  complete  log-likelihood,
conditional on the observed sequences, under a current parameter value
$\theta'$, leads to
\begin{align*}
\espp(\log \pr(\bX, \bY,  L,\e_{1:L}) |\bX, \bY ) = &\sum_{s=n\vee
  m}^{n+m} \prp(L=s |\bX, \bY )
 \Big\{ \sum_{u\in \E} \prp(\e_1=u |\Fs)\log \pi_{u} \\
&+\sum_{k=2}^{s}\sum_{u,v\in \E^2}
\prp((\e_{k-1},\e_{k})=(u,v) |\Fs) \log\pi_{uv} \\
&+\sum_{k=1}^s\sum_{a\in \A} \big[
\prp((\e, X)_k=(I_X,a) |\Fs) \log f(a) 
 +\prp((\e,Y)_k=(I_Y,a)|\Fs) \log g(a) \big]\\
& + \sum_{k=1}^s
\sum_{a,b} \prp(\e_{k-1}\neq M,
(\e, X, Y)_k=(M,a,b) |\Fs)%\\ & \qquad \times 
\log h(a,b) \\
& + \sum_{k=2}^s \sum_{a,b,c,d} 
 \prp( (\e_k,\e_{k-1} ,(X, Y)_k, (X, Y)_{k-1})=
(M,M,a,b,c,d) |\Fs)  \log \tilde{h}(a,b|c,d) 
\Big\}.
\end{align*}
% \begin{align*}
% &\espp(\log \pr(\bX, \bY,  L,\e_{1:L}) |\bX, \bY ) = \sum_{s=n\vee
%   m}^{n+m} \prp(L=s |\bX, \bY )\\
% &\times \Big\{ \sum_{u\in \E} \prp(\e_1=u |\Fs)\log \pi_{u} \\
% &+\sum_{k=2}^{s}\sum_{u,v\in \E^2}
% \prp((\e_{k-1},\e_{k})=(u,v) |\Fs) \log\pi_{uv} \\
% &+\sum_{k=1}^s\sum_{a\in \A} \big[
% \prp((\e, X)_k=(I_X,a) |\Fs) \log f(a) \\
% &\qquad +\prp((\e,Y)_k=(I_Y,a)|\Fs) \log g(a) \big]\\
% & + \sum_{k=1}^s
% \sum_{a,b} \prp(\e_{k-1}\neq M,
% (\e, X, Y)_k=(M,a,b) |\Fs)%\\ & \qquad \times 
% \log h(a,b) \\
% & + \sum_{k=2}^s \sum_{a,b,c,d} \\
% & \prp( (\e_k,\e_{k-1} ,(X, Y)_k, (X, Y)_{k-1})=
% (M,M,a,b,c,d) |\Fs) \\
% &\qquad \times \log \tilde{h}(a,b|c,d) 
% \Big\}.
% \end{align*}
There  are  two  main  issues  at  stake  here.  First,  the  quantity
$\prp(L_{n m}|X_{1:n},Y_{1:m})$  is not  given by  the forward-backward
equations. Second,  computing this conditional  expectation would also
require       the        knowledge       of       the       quantities
$$\prp(\e_{k-1},\e_{k} |\bX, \bY,L=s).$$  However, the forward
backward equations rather give access to quantities as for instance
\begin{equation*}
  \alpha^{I_X}(i-1,j-1)          \pi_{I_XM}         h(X_{i},Y_j)\beta^M(i,j)
  \\=\pr(\e_{k-1}=I_X,\e_{k}=M|\bX, \bY, Z_{k}=(i,j)).
\end{equation*}
In other  words, the conditional distribution  of $(\e_{k-1},\e_{k})$ is
known only conditional on the extra knowledge of the position $Z_k$ of the
path in the  lattice at time $k$. In conclusion,  in pair-HMM,
the  forward-backward  equations are  not  sufficient  to compute  the
expectation  of  the   complete  log-likelihood,  conditional  on  the
observed sequences.

\label{lastpage}

\end{document}